\documentclass[tikz,border=2pt,12pt]{amsart}

\usepackage{amsmath,amsthm,amssymb,mathrsfs,amsfonts,verbatim,enumitem,color,leftidx}
\usepackage{mathabx}
\usepackage{etoolbox,tikz} 
\usepackage{bbm}
\usepackage[all,tips]{xy}
\usepackage{graphicx,ifpdf}
\usepackage{stmaryrd,here}
\usetikzlibrary{positioning, fit, calc}
\ifpdf
\fi
\usepackage{hyperref}
\usepackage[right,displaymath,mathlines]{lineno}

\newtheorem{thm}{Theorem}[section]

\newtheorem{prop}[thm]{Proposition}

\theoremstyle{definition}
\newtheorem{de}[thm]{Definition}

\theoremstyle{remark}

\numberwithin{equation}{section}

\makeatletter

\newcommand{\Rmnum}[1]{\expandafter\@slowromancap\romannumeral #1@}
\makeatother

\begin{document}

\title[Critical exponents]{On the set of critical exponents of discrete groups acting on regular trees}


\author{Sanghoon Kwon}
\email{shkwon1988@gmail.com}

\thanks{}



\date{}



\begin{abstract}
We study the set of critical exponents of discrete groups acting on regular trees. We prove that for every real number $\delta$ between $0$ and $\frac{1}{2}\log q$, there is a discrete subgroup $\Gamma$ acting without inversion on a $(q+1)$-regular tree whose critical exponent is equal to $\delta$. Explicit construction of edge-indexed graphs corresponding to a quotient graph of groups are given.
\end{abstract}

\maketitle
\tableofcontents
\section{Introduction}

Let $\mathcal{T}_{q+1}$ be a $(q+1)$-regular tree and let $\textrm{Aut}(\mathcal{T}_{q+1})$ be the group of isometries of $\mathcal{T}_{q+1}$.  Denote by $V\mathcal{T}_{q+1}$ the set of all the vertices of $\mathcal{T}_{q+1}$. For a discrete subgroup $\Gamma$ of $\textrm{Aut}(\mathcal{T}_{q+1})$, the \emph{critical exponent} of $\Gamma$ is defined by
$$\delta_\Gamma=\underset{n\to\infty}{\overline{\lim}}\frac{\log\#\{\gamma\in\Gamma\colon d(v,\gamma v)\le n\}}{n}$$ for a fixed vertex $v\in V\mathcal{T}_{q+1}$. This does not depend on the choice of $v\in V\mathcal{T}_{q+1}$.

Equivalently, it can be defined as the infimum of $s>0$ such that the \emph{Poincar\'e} \emph{series} $\sum_{\gamma\in\Gamma}e^{-sd(v,\gamma v)}$ converges.

In some cases, it captures the dimension at ``infinity" of $\Gamma$ as well as measures the exponential growth rate of the size of $\Gamma$-orbits. It is shown in (\cite{Pa}) that if $\Gamma$ acts on $\mathcal{T}_{q+1}$ without a global fixed point on $\partial_\infty\mathcal{T}_{q+1}$, then $\delta_\Gamma$ is equal to the Hausdorff dimension of the \emph{conical limit set}. Here, $\partial_\infty\mathcal{T}_{q+1}$ is the \emph{Gromov boundary at infinity} of $\mathcal{T}_{q+1}$ defined as the set of equivalence class $[\ell]$ of geodesic ray $\ell$ starting from a fixed vertex $v\in V\mathcal{T}_{q+1}$, where two geodesic rays $\ell$ and $\ell'$ are equivalent if and only if $\{d_{\mathcal{T}}(\ell(n),\ell'(n))\colon n\in\mathbb{Z}_{>0}\}$ is bounded above. A conical limit set is the set of points $\omega$ in $\partial_\infty\mathcal{T}_{q+1}$ such that there is a sequence of points in $\Gamma v$ converging towards $\omega$ while staying at bounded distance from a geodesic ray ending at $\omega$.


If the quotient graph $\Gamma\backslash\mathcal{T}_{q+1}$ is finite, then we say that $\Gamma$ is a \emph{uniform lattice} in $\textrm{Aut}(\mathcal{T}_{q+1})$ and in this case $\delta_\Gamma=\log q$. If $\Gamma$ is either itself finite or has a cyclic group as a finite index subgroup, then $\delta_\Gamma$ is equal to zero. In general, we always have $0\le\delta_\Gamma\le\log q$. Recently, it is proved in \cite{CDS} that if $\Gamma$ is convex co-compact (that is, every limit point of $\Gamma$ is a conical limit point) acting isometrically on a $CAT(-1)$ space and $\Lambda$ is a subgroup of $\Gamma$, then $\delta_\Gamma=\delta_{\Lambda}$ if and only if $\Lambda$ is \emph{co-amenable} in $\Gamma$. Recall that $\Lambda$ is said to be co-amenable in $\Gamma$ if the left coset space $\Lambda\backslash\Gamma$ has $\Gamma$-invariant mean. See also \cite{No} and \cite{No2} for the similar results in regular and arbitrary graphs.

In this article, we ask the following question of existence: For a given number $\delta\in [0,\log q]$, is there a discrete subgroup $\Gamma$ of $\textrm{Aut}(\mathcal{T}_{q+1})$ such that $\delta_\Gamma=\delta$? We give a partial answer to this question. Namely, for $\delta\in [0,\frac{1}{2}\log q]$, there exists a discrete subgroup $\Gamma$ of $\textrm{Aut}(\mathcal{T}_{q+1})$ for which $\delta_\Gamma=\delta$.

Let $E=\{\delta\in[0,\log q]\colon \delta=\delta_\Gamma\textrm{ for some discrete }\Gamma<\textrm{Aut}(\mathcal{T}_{q+1})\}$.

\begin{thm}\label{main1}
If $\delta\in[0,\frac{1}{2}\log q]$, then $\delta\in E$.
\end{thm}

In Section 2, we summarize the connection between graph of groups and its edge-indexed graphs and their covering theory, which is the main idea of the proof. The proof of Theorem~\ref{main1} is given in Section 3 in which we construct an edge-indexed graph corresponding to a real number $\delta$ in $[0,\frac{1}{2}\log q]$, for which the fundamental group of its finite grouping has the critical exponent $\delta$. In Section 4, we give some remarks involving two questions: Which real numbers larger than $\frac{1}{2}\log q$ are realized? Which algebraic integers are realized as a critical exponents of free groups?

\subsubsection*{Acknowledgement}


\section{Edge-indexed graphs and graph of groups}

We briefly review the essential features of the covering theory for graph of groups which is the main idea of the proof. Being developed by Serre \cite{Se} and the substantial contribution of the subsequent work of Bass \cite{Ba}, the theory is currently called as Bass-Serre theory. We mainly follows \cite{Se} and refer to \cite{Ba} for further details.

Let $A$ be an undirected graph which is allowed to have loops and multiple edges. We denote by $VA$ the set of vertices of $A$ and by $EA$ the set of edges of $A$. Given an undirected graph $A$, we get a \emph{symmetric} directed graph corresponding to $A$, whose set of vertices are the same as $VA$ and every edge of $A$ are bidirected. Denote by $E^{ori}\!\!A$ the set of all oriented edges of the symmetric directed graph corresponding to $A$, hence $|E^{ori}A|=2|EA|$.

For $e\in E^{ori}\!\!A$, let $\overline{e}\in E^{ori}\!\!A$ be the opposite edge of $e$ and let $\partial_0e$ and $\partial_1e$ be the initial vertex and the terminal vertex of $e$, respectively. 

\begin{de}Let $i\colon E^{ori}\!\!A\to \mathbb{Z}_{>0}$ be a map assigning a positive integer to each oriented edge. We say $(A,i)$ an \emph{edge-indexed graph}.
\end{de}

\begin{de}By a \emph{graph of groups} $\mathbf{A}=(A,\mathcal{A})$, we mean a connected graph $A$ together with attached groups $\mathcal{A}_a$ $(a\in VA)$, $\mathcal{A}_e=\mathcal{A}_{\overline{e}}$ $(e\in E^{ori}\!\!A)$, and monomorphisms $\alpha_e\colon \mathcal{A}_e\to\mathcal{A}_{\partial_1 e}$ $(e\in E^{ori}\!\!A)$. 
\end{de}

An \emph{isomorphism} between two graph of groups $\mathbf{A}=(A,\mathcal{A})$ and $\mathbf{A'}=(A',\mathcal{A'})$ is an isomomorphism $\phi\colon A\to A'$ between two underlying graphs together with the set of isomomorphisms $\phi_a\colon \mathcal{A}_a\to \mathcal{A'}_{\phi(a)}$ and $\phi_e\colon \mathcal{A}_e\to \mathcal{A'}_{\phi(e)}$ satisfying the following property: for each $e\in EA$, there is an element $h_e\in \mathcal{A'}_{\phi(\partial_1e))}$ such that 
$$\phi_{\partial_1e}(\alpha_e(g))=h_e\cdot(\alpha'_{\phi(e)}(\phi_e(g)))\cdot h_e^{-1}$$
for all $g\in\mathcal{A}_e$ (\cite{Ba}, Section 2).

Given an edge-indexed graph $(A,i)$, a graph of groups $(A,\mathcal{A})$ is called a \emph{grouping} of $(A,i)$ if $i(e)=[\mathcal{A}_{\partial_1e}\colon\alpha_e\mathcal{A}_e]$ and called a \emph{finite grouping} if all $\mathcal{A}_a$ $(a\in VA)$ are finite.

Suppose that we have a graph of groups $\mathbf{A}$. Choosing a base point $a_0\in VA$, we can define a \emph{fundamental group} $\pi_1(\mathbf{A},a_0)$ (\cite{Se} Section 5.1), a \emph{universal covering tree} $(\widetilde{\mathbf{A},a_0})$ and an action without inversion of $\pi_1(\mathbf{A},a_0)$ on $(\widetilde{\mathbf{A},a_0})$ with a morphism $p\colon (\widetilde{\mathbf{A},a_0}) \to A$ which can be identified with the quotient projection (\cite{Se}, Section 5.3).

\begin{de} Given a graph of groups $\mathbf{A}=(A,\mathcal{A})$, we denote by $F(\mathcal{A}, E^{ori})$ the group generated by the groups $\mathcal{A}_a$, $(a\in VA)$ and the elements $e\in E^{ori}$, subject to the relations
$$\overline{e}=e^{-1}\textrm{ and }e\alpha_e(g)e^{-1}=\alpha_{\overline{e}}(g)\textrm{ for }e\in E^{ori}A\textrm{ and }a\in \mathcal{A}_e.$$
Let $\pi_1(\mathbf{A},a_0)$ be the set of elements of $F(\mathcal{A},E^{ori})$ of the form $$g_0e_1g_1e_2\cdots e_{n-1}g_{n-1},$$ where $o(e_{i+1})=t(e_i)$ (mod $n)$ and $g_i\in o(e_{i+1})$. It is a subgroup of $F(\mathcal{A},E^{ori})$, called \emph{fundamental group} of $\mathbf{A}$ based at $a_0$.
\end{de}

Given a graph of groups $\mathbf{A}$, the graph $\widetilde{X}=(\widetilde{\mathbf{A},a_0})$ is defined as
$$V\widetilde{X}=\bigcup_{a\in VA}\pi_1(\mathbf{A},a_0)/\textrm{Stab}_{\pi_1(\mathbf{A},a_0)}(a)$$ and $$E\widetilde{X}=\bigcup_{e\in E^{ori}A}\pi_1(\mathbf{A},a_0)/\textrm{Stab}_{\pi_1(\mathbf{A},a_0)}(e).$$

\begin{thm}[\cite{Se}, Theorem 12] The graph $\widetilde{X}$ defined as above is a tree.
\end{thm}
The tree $\widetilde{X}=(\widetilde{\mathbf{A},a_0})$ is called a \emph{universal covering tree} of the graph of groups $\mathbf{A}$.

Suppose that a group $\Gamma$ acts on a graph $A$. We call $s\in\Gamma$ an \emph{inversion} on $A$ if $se=\overline{e}$ for some $e\in E^{ori}\!\!A$. If $\Gamma$ acts without inversions, then the quotient graph $\Gamma\backslash A$ is well-defined and we have a natural projection $p\colon A\to\Gamma\backslash A$.

When a group $\Gamma$ acts without inversion on a tree $\mathcal{T}_{q+1}$, we can naturally identify $\Gamma$ and $\mathcal{T}_{q+1}$ with the fundamental group and universal covering tree, respectively, of a graph of groups called \emph{quotient graph of groups}.

\begin{de}
Suppose that $\Gamma$ acts without inversion on a tree $\mathcal{T}_{q+1}$. For each $v\in V(\Gamma\backslash\mathcal{T})$ and $e\in E^{ori}(\Gamma\backslash\mathcal{T})$, we choose any corresponding vertex $\widetilde{v}\in V\mathcal{T}$ and edge $\widetilde{e}\in E\mathcal{T}$. Let $\overline{\widetilde{e}}=\widetilde{\overline{e}}$ and fix an element $\gamma_e\in\Gamma$ which satisfies $\gamma_e\widetilde{\partial_1(e)}=\partial_1(\widetilde{e})$. Define $\mathcal{A}_v$ and $\mathcal{A}_e$ be the stabilizer of $\widetilde{v}$ and $\widetilde{e}$ in $\Gamma$, respectively, and let $\alpha_e\colon\mathcal{A}_e\to\mathcal{A}_{\partial_1e}$ be the monomorphism $h\mapsto \gamma_e^{-1}h\gamma_e$. Then the graph of groups $(\Gamma\backslash\mathcal{T}_{q+1},\mathcal{A})$ is called the \emph{quotient graph of groups}.
\end{de}

It does not depend on the choice of $\widetilde{v},\widetilde{e}$ and $\gamma_e$ up to isomorphism of graph of groups (\cite{Ba}, Section 3).

Conversely, if $(A,\mathcal{A})$ is a graph of groups and $(A,i)$ is the corresponding edge-indexed graph, then fixing a basepoint $a_0\in VA$, the universal covering tree $\mathcal{T}=\widetilde{(A,a_0)}$ and the natural projection $\pi\colon \mathcal{T}\to A$ depend only on the edge indexed graph $(A,i)$ (\cite{Ba}, Remark 1.18).

\section{Proof of Theorem}

Our goal is to construct an edge-indexed graph $(A,i)$ which admits a finite grouping $\mathbf{A}$ of which the fundamental group has critical exponent $\delta\in [0,\frac{1}{2}\log q]$.

We recall the following criterion for an edge-indexed graph to have a finite grouping.

\begin{thm}[\cite{bk90}, Corollary 2.4] An edge-indexed graph $(A,i)$ admits a finite grouping if and only if there is a positive integral valued function $N\colon VA\to\mathbb{Z}_{>0}$ on $VA$ satisfying \begin{equation}\label{ivo}\frac{i(e)}{i(\overline{e})}=\frac{N(\partial_0e)}{N(\partial_1e).}\end{equation}
\end{thm}

We say such a function $N$ an \emph{integral vertex ordering}.

 For convenience to describe some graphs constructed in the proof, we introduce a graph operation between two graphs which we call \emph{graph merging}.
\begin{de} Let $G_1=(V_1,E_1)$ and $G_2=(V_2,E_2)$ be two (undirected or directed) graphs. The \emph{graph merging} $G_1\underset{x,y}{\star} G_2$ with respect to $x\in V_1$ and $y\in V_2$ is a graph $G=(V_1\cup V_2/{x\sim y}, E_1\cup E_2)$, obtained by identifying two vertices $x$ and $y$ in $G_1\cup G_2$.
\end{de}

\vspace{.1 in}
\begin{figure}[H]
\begin{center}
\begin{tikzpicture}[every loop/.style={}]
  \tikzstyle{every node}=[inner sep=0pt]
  \node (0) {$\bullet$} node [left=22pt] at (0,0) {$G_1$}; 
  \node (0) {}node [above=4pt] at (0,0) {$x\sim y$};
  \node (1) {}node [right=24pt] at (0,0) {$G_2$};
  \draw[dashed] (-1cm,0cm) circle (1cm) ;
  \draw[dashed] (1cm,0cm) circle (1cm) ;
\end{tikzpicture}
\vspace{0.5em}
\caption{The graph $G_1\underset{x,y}{\star} G_2$}
\end{center}
\end{figure}
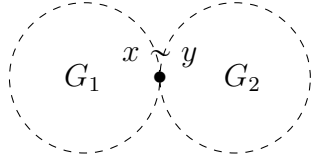

We give a proof of Theorem~1.1. Let $A_0$ be a ray given as follows:

\begin{figure}[H]
\begin{center}
\begin{tikzpicture}[every loop/.style={}]
  \tikzstyle{every node}=[inner sep=0pt]
  \node (0) {$\bullet$} node [below=4pt] at (0,0) {$x_0$};
  \node (2) at (1.5,0) {$\bullet$} node [below=4pt] at (1.5,0) {$x_1$}; 
  \node (4) at (3,0) {$\bullet$}node [below=4pt] at (3,0) {$x_2$}; 
  \node (6) at (4.5,0) {$\bullet$}node [below=4pt] at (4.5,0) {$x_3$}; 
  \node (8) at (6,0) {$\bullet$}node [below=4pt] at (6,0) {$x_4$}; 
  \node (10) at (7.5,0) {$\bullet$}node [below=4pt] at (7.5,0) {$x_5$}; 
  \node (11) at (9.3,0) {$\bullet \cdots$}node [below=4pt] at (9,0) {$x_6$}; 

  \path[-] (0) edge node [above=4pt] {} (2)
 (2) edge node [above=4pt] {} (4)
 (4) edge node [above=4pt] {} (6)
 (6) edge node [above=4pt] {} (8)
 (8) edge node [above=3.2pt] {} (10)
 (10) edge node [above=4pt] {} (11);
\end{tikzpicture}
\vspace{0.5em}
\caption{A ray $A_0$}
\end{center}
\end{figure}
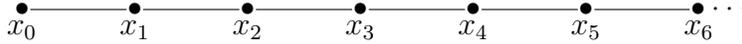Let $I\cup J=\mathbb{Z}_{>0}$ be an arbitrary partition of $\mathbb{Z}_{>0}$. If $n\in I$, then we assign $i(\overrightarrow{x_nx_{n-1}})=q$ and $i(e)=1$ otherwise. We adjoin a tree $\mathcal{T}_{q+1}'$ which is $(q+1)$-regular except one $y$ vertex of degree $q$ to $x_0$.  If $n\in J$, then for each $n$ we adjoin a copy of a tree $\mathcal{T}_{q+1}''$ which is $(q+1)$-regular except one vertex $z_n$ of degree $q-1$ to $x_n$. In other words, $(A,i)=(A_0\underset{x_0,y}{\star}\mathcal{T}_{q+1}'\underset{x_j,z_j,j\in J}{\star}\mathcal{T}_{q+1}'',i)$, where $i(\overrightarrow{x_nx_{n-1}})=q$ for each $n\in I$ and $i(e)=1$ otherwise.

Then $(A,i)$ admits an integral vertex ordering. Indeed, we can define $N(x_0)=1$ and $N(x_j)=i(\overrightarrow{x_jx_{j-1}})N(x_{j-1})$ for all $j\ge 1$. Then $N$ satisfies the equation~(\ref{ivo}).

Let $\Gamma$ be the subgroup of $\textrm{Aut}(\mathcal{T}_{q+1})$ such that the edge-indexed graph of the quotient graph of groups $(\Gamma\backslash \mathcal{T}_{q+1},\mathcal{A})$ is $(A,i)$. As we discussed in Section 2, such a group $\Gamma$ exists and is isomorphic to the fundamental group of a graph of groups obtained by certain finite grouping of $(A,i)$. In addition, for any finite grouping, the group $\Gamma$ is discrete since $\Gamma_x$ is finite for every $x\in VA$.

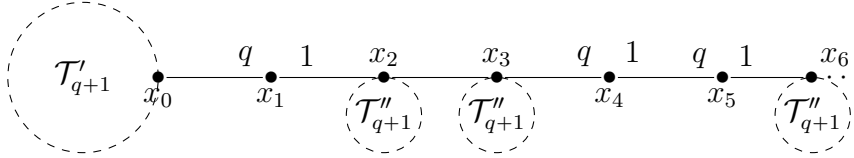
\begin{figure}[H]
\begin{center}
\begin{tikzpicture}[every loop/.style={}]
  \tikzstyle{every node}=[inner sep=0pt]
  \node (0) {$\bullet$} node [below=4pt] at (0,0) {$x_0$};
  \node (2) at (1.5,0) {$\bullet$} node [below=4pt] at (1.5,0) {$x_1$}; 
  \node (4) at (3,0) {$\bullet$}node [above=4pt] at (3,0) {$x_2$}; 
  \node (6) at (4.5,0) {$\bullet$}node [above=4pt] at (4.5,0) {$x_3$}; 
  \node (8) at (6,0) {$\bullet$}node [below=4pt] at (6,0) {$x_4$}; 
  \node (10) at (7.5,0) {$\bullet$}node [below=4pt] at (7.5,0) {$x_5$}; 
  \node (11) at (9,0) {$\bullet \cdots$}node [above=4pt] at (9,0) {$x_6$}; 
  \node (-1) at (-1,0) {$\mathcal{T}_{q+1}'$};
  \node (12) at (3,-0.5) {$\mathcal{T}_{q+1}''$};
  \node (13) at (4.5,-0.5) {$\mathcal{T}_{q+1}''$};
  \node (14) at (8.7, -0.5) {$\mathcal{T}_{q+1}''$};
  \draw[dashed] (-1cm,0cm) circle (1cm);
  \draw[dashed] (3cm,-0.5cm) circle (0.5cm);
  \draw[dashed] (4.5cm,-0.5cm) circle (0.5cm);
  \draw[dashed] (8.7cm,-0.5cm) circle (0.5cm);

  \path[-] (0) edge node [above=4pt] {\quad \quad $q$} (2)
 (2) edge node [above=4pt] {$1$ \quad \quad } (4)
 (4) edge node [above=4pt] {} (6)
 (6) edge node [above=4pt] {\quad\quad $q$} (8)
 (8) edge node [above=3.2pt] {$1$ \quad \,\,$q$} (10)
 (10) edge node [above=4pt] {$1$ \quad \quad} (11);
\end{tikzpicture}
\vspace{0.5em}
\caption{$(A,i)$, Example: $1,4,5\in I$ and $2,3,6\in J$}
\end{center}
\end{figure}

For each $n\in\mathbb{Z}_{>0}$, if we denote by $\Gamma_{x}$ the stabilizer of $x$ in $\Gamma$ and by $\widetilde{x_0}$ a lift of $x_0$ in the universal covering tree $\mathcal{T}_{q+1}$, then 
 \[
  \#\{\gamma\in\Gamma\colon d(\widetilde{x_0},\gamma \widetilde{x_0}) =2n\}= \left\{\begin{array}{cl}
          \displaystyle 0& \textrm{if }N(x_n)=N(x_{n-1})\\
        \displaystyle \frac{q-1}{q}N(x_n)|\Gamma_{x_0}|& \textrm{if }N(x_n)=qN(x_{n-1})
        \end{array}\right.
  .\]

Let $s_n=\#\{1\le k\le n\colon i(\overrightarrow{x_kx_{k-1}})=q\}$.
Then, \begin{align*}\displaystyle\delta_\Gamma=\underset{n\to\infty}{\overline{\lim}}\frac{1}{2n}\log\frac{(q-1)q^{s_n}}{q}=\underset{n\to\infty}{\overline{\lim}}\frac{s_n}{2n}\log q.
\end{align*}
By taking a suitable partition $I$ and $J$, we can make $\displaystyle\underset{n\to\infty}{\overline{\lim}}\frac{s_n}{n}$ be any real number in $[0,1]$. This completes the proof of Theorem~\ref{main1}.

\section{Further discussion: Ihara zeta function and generating functions}

We give some remarks towards two problems in this section. First question is to characterize explicitly the set $$E=\{\delta\in[0,\log q]\colon \delta=\delta_\Gamma\textrm{ for some discrete }\Gamma<\textrm{Aut}(\mathcal{T}_{q+1})\}.$$

We are interested in constructing a discrete group for a real number $\delta\in\left(\frac{1}{2}\log q,\log q\right)$ whose critical exponent is eqaul to $\delta$. We consider an edge-indexed graph obtained by the graph merging of those appeared in Section 3. The critical exponent of the discrete group corresponding to the one is related to those of smaller edge-indexed graphs, being the radius convergence of certain Laurent series.

Given an edge-indexed graph $X=(A,i)$, let $N_X(m)$ be the number of closed paths of length $m$ without backtracking or tails in the graph $A$ with weights $i$. More precisely, given a closed cycle $C=e_1\cdots e_m$ of length $m$, let $w(C)=i(e_1)+\cdots i(e_m)$ and
$$N_X(m)=\sum_{C=e_1\cdots c_m}w(C)$$
where the summation runs over the closed cycle of length $m$ without backtracking or tails.

Let us denote by $F_X(u)$ the generating function given by $$\displaystyle \sum_{m=1}^{\infty} N_X(m)u^m.$$
Suppose that every closed cycle of $X$ and $Y$ passes through the vertex $x$ and $y$, respectively.

Then, we have
\begin{align*}N_{X\underset{x,y}{\star}Y}(m)=&\left(N_X(m)+\sum N_{X}(k_1)N_{Y}(k_2)+\sum N_{X}(k_1)N_{Y}(k_2)N_{X}(k_3)+\cdots\right)\\+&\left(N_Y(m)+\sum N_{Y}(k_1)N_{X}(k_2)+\sum N_{Y}(k_1)N_{X}(k_2)N_{Y}(k_3)+\cdots\right) 
\end{align*}
where all the summations run over $k_1+k_2+\cdots k_n=m$. Hence, it follows that
\begin{align*}&\sum_{m=1}^{\infty}N_{X\underset{x,y}{\star}Y}(m)u^m\\=&\left(\sum_{m=1}^{\infty}N_X(m)u^m+\sum_{m=1}^{\infty}N_X(m)u^m\sum_{m=1}^{\infty}N_Y(m)u^m+\cdots\right)\\
+&\left(\sum_{m=1}^{\infty}N_Y(m)u^m+\sum_{m=1}^{\infty}N_Y(m)u^m\sum_{m=1}^{\infty}N_X(m)u^m+\cdots\right)
\end{align*}
which yields the following proposition. 

\begin{prop} Let $X$ and $Y$ be an edge-indexed graph such that every closed cycle of $X$ and $Y$ passes through the vertex $x$ and $y$, respectively. Then we have
$$F_{X\star Y}(u)=\sum_{m=1}^{\infty}N_{X\underset{x,y}{\star}Y}(m)u^m=\frac{F_X(u)+F_Y(u)+2F_X(u)F_Y(u)}{1-F_X(u)F_Y(u)}.$$
\end{prop}
For instance, if both $X$ and $Y$ are edge-indexed graphs given in Section 3, then $N_X(m)$ is either $q^k$ for some $k<m$ or is equal to $0$. Thus, given any $0$-$1$ sequences $a_n$ and $b_n$, if we denote by $s_n$ and $t_n$ the partial sum of each, then any real numbers realized as a solution of the equation $u\in\mathbb{R}$ such that $$\sum_{m=1}^{\infty} q^{s_m}u^m\sum_{m=1}^{\infty} q^{t_m}u^m=1$$ are contained in the set $E$.

\vspace{.1 in}
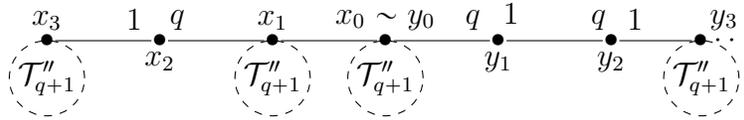
\begin{figure}[H]
\begin{center}
\begin{tikzpicture}[every loop/.style={}]
  \tikzstyle{every node}=[inner sep=0pt]
  \node (0) {$\bullet$} node [above=4pt] at (0,0) {$x_3$};
  \node (2) at (1.5,0) {$\bullet$} node [below=4pt] at (1.5,0) {$x_2$}; 
  \node (4) at (3,0) {$\bullet$}node [above=4pt] at (3,0) {$x_1$}; 
  \node (6) at (4.5,0) {$\bullet$}node [above=4pt] at (4.5,0) {$x_0\sim y_0$}; 
  \node (8) at (6,0) {$\bullet$}node [below=4pt] at (6,0) {$y_1$}; 
  \node (10) at (7.5,0) {$\bullet$}node [below=4pt] at (7.5,0) {$y_2$}; 
  \node (11) at (9,0) {$\bullet \cdots$}node [above=4pt] at (9,0) {$y_3$}; 
  \node (-1) at (0,-0.5) {$\mathcal{T}_{q+1}''$};
  \node (12) at (3,-0.5) {$\mathcal{T}_{q+1}''$};
  \node (13) at (4.5,-0.5) {$\mathcal{T}_{q+1}''$};
  \node (14) at (8.7, -0.5) {$\mathcal{T}_{q+1}''$};
  \draw[dashed] (0cm,-0.5cm) circle (0.5cm);
  \draw[dashed] (3cm,-0.5cm) circle (0.5cm);
  \draw[dashed] (4.5cm,-0.5cm) circle (0.5cm);
  \draw[dashed] (8.7cm,-0.5cm) circle (0.5cm);

  \path[-] (0) edge node [above=4pt] {\quad \quad $1$} (2)
 (2) edge node [above=4pt] {$q$ \quad \quad \, } (4)
 (4) edge node [above=4pt] {} (6)
 (6) edge node [above=4pt] {\quad\quad $q$} (8)
 (8) edge node [above=4pt] {$1$ \quad\quad $q$} (10)
 (10) edge node [above=4pt] {$1$ \quad \quad} (11);
\end{tikzpicture}
\vspace{0.5em}
\caption{$X\underset{x_0,y_0}{\star} Y\underset{x_0}{\star}\mathcal{T}_{q+1}''$}
\end{center}
\end{figure}

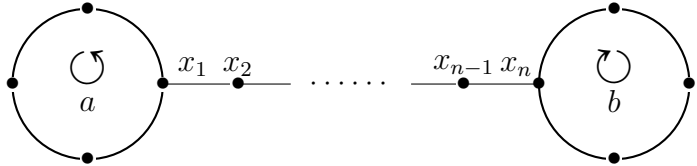
\begin{figure}[htbp]
\begin{center}
\begin{tikzpicture}[every loop/.style={}]
  \tikzstyle{every node}=[inner sep=0pt]
  \draw[thick] (0cm,0cm) circle (1cm);
  \draw (-1,0) node[fill=white]{$\bullet$};
  \draw (0,1) node[fill=white]{$\bullet$};
  \draw (0,-1) node[fill=white]{$\bullet$};
  \draw (0,0) node{\Large{$\underset{a}{\circlearrowleft}$}};
  \draw (1,0) node[fill=white]{$\bullet$}; 
  \draw (1cm,0cm) --(2cm,0cm);
  \draw (2,0) node[fill=white]{$\bullet$};
  \draw (2,0)--(2.7,0);
  \draw (3.5,0) node{$\cdots\cdots$};
  \draw (4.3,0) --(5cm,0cm); 
  \draw (1.4cm,0cm) node[above=3pt]{$x_1$};
  \draw (2cm,0cm) node[above=3pt]{$x_2$};
  \draw (5cm,0cm) node[above=3pt]{$x_{n-1}$};
  \draw (5.7cm,0cm) node[above=3pt]{$x_{n}$};
  \draw (5,0) node[fill=white]{$\bullet$}; 
  \draw (6,0) node[fill=white]{$\bullet$};
  \draw (5,0)--(6,0);
  \draw[thick] (7cm,0cm) circle (1cm);
  \draw (7,1) node[fill=white]{$\bullet$};
  \draw (8,0) node[fill=white]{$\bullet$};
  \draw (7,-1) node[fill=white]{$\bullet$};
  \draw (7,0) node{\Large{$\underset{b}{\circlearrowright}$}};
\end{tikzpicture}\label{dumbbell}
\caption{Dumbbell graph $D_{a,b,n}$}
\end{center}
\end{figure}

The other question concerns the set of critical exponents of \emph{free groups}. Using Ihara zeta function (the definition will be given below), we see that if $\Gamma$ is a free subgroup of $\textrm{Aut}(\mathcal{T}_{q+1})$, then the critical exponent $\delta_\Gamma$ is equal to $\log \alpha$ for some algebraic integer $\alpha$. We are interested in characterizing the irreducible polynomials of such $\alpha$. For example, in \cite{Kw} the critical exponents of certain dumbbell graphs, consisting of two vertex-disjoint cycles joining them having only its end-vertices in common with the two cycles (see Figure~\ref{dumbbell}), are described.

We briefly explain how to calculate the irreducible polynomial of the critical exponent of the free group $\Gamma$. This is related to the radius of convergence of the zeta function of the graph $\Gamma\backslash\mathcal{T}_{q+1}^\textrm{min}$. Here, $\mathcal{T}_{q+1}^\textrm{min}$ is the minimal $\Gamma$-invariant subtree of $\mathcal{T}_{q+1}$.

Suppose that $A$ is a finite connected undirected graph with no degree 1 vertices. A closed cycle is called \emph{primitive} if it is not a finite multiple of a smaller closed cycle. Let $P=(e_1,e_2,\ldots,e_{l(P)-1},e_{l(P)})$ be a primitive closed cycle without backtracking. That is, $i(e_1)=t(e_{l(P)})$, $e_{i+1}\ne e_i^{-1}$ $(\textrm{mod }l(P))$ for all $i$ and $P\ne D^m$ for any integer $m\ge 2$ and a closed cycle $D$ in $A$. If a closed cycle $Q$ is obtained by changing the cyclic order of $P$, then we say $P$ and $Q$ are \emph{equivalent}. A \emph{prime} $[P]$ in $A$ are equivalence classes of primitive closed cycle without backtracking. 

\begin{de} The \emph{Ihara zeta function} of a finite graph $A$ is defined at $u\in\mathbb{C}$, for which $|u|$ is sufficiently small, by
$$Z_A(u)=\prod_{[P]}(1-u^{l(P)})^{-1}$$ where $[P]$ runs over the primes of $A$.
\end{de}

The Ihara determinant formula (\cite{Ih} which was established also by \cite{Ba2} and \cite{KS}) says that 
$$Z_A(u)=\frac{1}{(1-u^2)^{\chi(A)-1}\det(I-Adj(u)+Qu^2)}$$
where $\chi=|EA|-|VA|+1$, $Adj$ is the vertex adjacency matrix of $A$ and $Q$ is diagonal matrix whose $j$-th diagonal entry is $\textrm{deg}(v_j)-1$.

Let $\Delta_A$ be the g.c.d. of $\{l(P)\colon P\textrm{ prime in }A\}$ and $\pi_A(n)$ be $$\#\{[P]\colon P\textrm{ prime in }A, l(P)=n\}.$$ The following is the graph version of the prime geodesic theorem.

\begin{thm}[\cite{Te}, Theorem 10.1] Let $R_A$ be the radius of convergence of $Z_A(u)$. If $\Delta_A=1$, then $$\lim_{n\to\infty}n\pi(n)R_A^{n}=1.$$ If $\Delta_A>1$, then $\pi(m)=0$ unless $\Delta_A|m$ and $$\lim_{n\to\infty}n\Delta_A\pi(n\Delta_A) R_A^{n\Delta_A}=1.$$ 
\end{thm}

In the proof of the above theorem (\cite{Te}, Theorem 10.1), it is shown that if we denote by $W$ the adjacency matrix of the oriented line graph of $A$, then we have 
$$\pi_A(n)\sim\frac{1}{n}\sum_{\underset{\lambda\in \textrm{Spec}W}{|\lambda|\,\textrm{max}}}\lambda^n\quad\textrm{and}\quad N_A(m)=\sum_{\lambda\in\textrm{Spec}W}\lambda^m$$ which implies that $$\underset{n\to\infty}{\overline{\lim}}\frac{\log \pi_A(n)}{n}=\underset{n\to\infty}{\overline{\lim}}\frac{\log N_A(n)}{n}=\log\lambda_A^{\textrm{max}},$$
where $\lambda_A^{\textrm{max}}\in \textrm{Spec}\,W$ is the eiegenvalue of $W$ with maximum modulus. In \cite{KS}, we also have $$Z_A(u)=\frac{1}{\det(I-uW)}.$$

Suppose that $\Gamma$ is a free subgroup of $\textrm{Aut}(\mathcal{T}_{q+1})$. Then, we have 
\begin{align*}\delta_\Gamma=&\underset{n\to\infty}{\overline{\lim}}\frac{\log\#\{\gamma\in\Gamma\colon d(v,\gamma v)\le n\}}{n}\\=&\underset{n\to\infty}{\overline{\lim}}\frac{\log|\textrm{Stab}_\Gamma(\widetilde{v})|\sum_{k=1}^{n}N_A(k)}{n} \\=&\underset{n\to\infty}{\overline{\lim}}\frac{\log N_A(n)}{n}\\=&
 \underset{n\to\infty}{\overline{\lim}}\frac{\log \pi_A(n)}{n}=\lambda_A^{\textrm{max}}=\log R_A^{-1}
\end{align*}
where we take $A$ by $\Gamma\backslash\mathcal{T}_{q+1}^\textrm{min}$.

\end{document}